\documentclass[journal]{IEEEtran}

\usepackage{cite}
\usepackage{algorithmic}
\usepackage{array}
\usepackage{url}
\usepackage{booktabs}
\usepackage[acronym]{glossaries}

\usepackage{xcolor}
\definecolor{revision}{rgb}{0,0,1}

\usepackage{amsmath}
\usepackage{amssymb}
\usepackage{subcaption}

\ifCLASSINFOpdf
  \usepackage[pdftex]{graphicx}
\else
  \usepackage[dvips]{graphicx}
\fi

\usepackage{todonotes}

\newacronym[plural=AVs]{av}{AV}{automated vehicle}
\newacronym[plural=CVs]{cv}{CV}{connected vehicle}
\newacronym[plural=CAVs]{cav}{CAV}{connected and automated vehicle}
\newacronym[plural=SAVs]{sav}{SAV}{shared autonomous vehicle}
\newacronym[plural=SAEVs]{saev}{SAEV}{shared autonomous electric vehicle}

\begin{document}
%
\title{Enhanced Mobility with Connectivity and Automation: A Review of Shared Autonomous Vehicle Systems}
%
%
%

\author{Liuhui~Zhao,~\IEEEmembership{Member,~IEEE,}
        Andreas~A.~Malikopoulos,~\IEEEmembership{Senior~Member,~IEEE}
\thanks{This research was supported in part by ARPAE's NEXTCAR program under the
award number DE-AR0000796 and by the Delaware Energy Institute (DEI).}%
\thanks{The authors are with Department
of Mechanical Engineering, University of Delaware, Newark, DE 30332; e-mails: (lhzhao@ieee.org, andreas@udel.edu).}
}

\markboth{IEEE Intelligent Transportation Systems Magazine}%
{}
%


\maketitle

\begin{abstract}
Shared mobility can provide access to transportation on {a custom basis} without vehicle ownership. The advent of connected and automated vehicle technologies can further enhance the potential benefits of shared mobility systems. Although the implications of a system with shared autonomous vehicles have been investigated, the research  reported in the literature has exhibited contradictory outcomes. In this paper, we present a summary of the research efforts in shared autonomous vehicle systems that have been reported in the literature to date and discuss potential future research directions.
\end{abstract}


%

\IEEEpeerreviewmaketitle
\section{Introduction}\label{sec:int}
\subsection{Movitation}
In a rapidly urbanizing world, we need to make fundamental transformations in how we use and access transportation. We are currently witnessing an increasing integration of our energy and transportation which, coupled with the human interactions, is giving rise to a new level of complexity \cite{Malikopoulos2016b} in emerging transportation systems such as \glspl{cav} and shared mobility. As we move to increasingly complex emerging transportation systems, new control approaches \cite{Malikopoulos2008a,Malikopoulos2009} are needed to optimize their impact on the mobility system behavior.

\begin{figure}[b]
    \centering
    \includegraphics[width=0.45\textwidth]{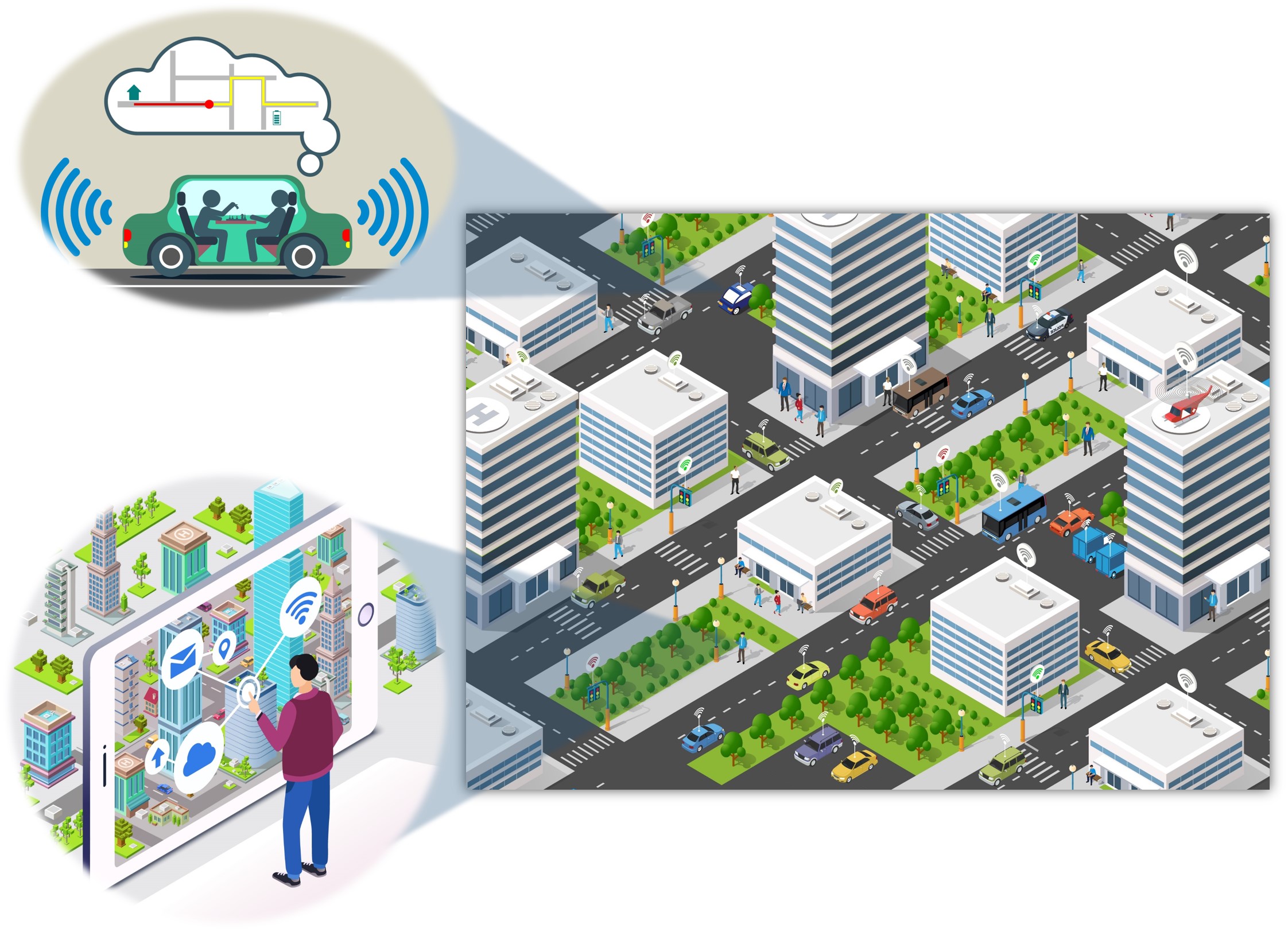}
    \caption{A view of a city enhanced by connectivity and automation.}
    \label{fig:shared_mobility}
\end{figure}

Shared mobility includes a variety of service models (e.g., carsharing, ridesharing, bikesharing) to meet travel needs and may result in a transformative impact on urban mobility \cite{Shaheen2016a, Shaheen2016b, SAE2018, cohen2018, standing2018} and landscape. {As shared mobility services evolve, there has been a debate on their potential impact \cite{feigon2016, cohen2018, shaheen2018}. The advent of intelligent transportation systems and information technologies has aimed at facilitating shared mobility services (Fig.~\ref{fig:shared_mobility}). In this context, impact analysis of the introduction of connected vehicles and \glspl{av} into existing shared mobility services is vital to identify the opportunities and challenges related to a shared autonomous mobility system.} In this paper, we review the research reported in the literature on carsharing enhanced by vehicle connectivity and automation technologies, i.e., \gls{sav} system, and discuss  potential implications in the environment and urban mobility.  

\subsection{Background}
{There are different types of carsharing service models, including round-trip carsharing, one-way station-based or free-floating carsharing, and peer-to-peer carsharing \cite{Shaheen2016b, Santos2016}.} In the past few years, short-term vehicle rental services provided by carsharing companies in major cities has attracted millions of users, while the number is expected to grow significantly \cite{feigon2016, shaheen2018, Martin2016}. Generally, it is believed that carsharing has positive impacts on energy use and greenhouse gas (GHG) emissions \cite{Martin2011a, Martin2011, Martin2016, firnkorn2015, Chen2016}, particularly when low-polluting vehicles are introduced into the transportation systems \cite{Barth2002}. {Although there is evidence that the use of carsharing services leads to a decrease in vehicle ownership \cite{Martin2011, Martin2011a, Martin2016}, location-specific variations (e.g., urban form, level of transit service, availability of alternative modes, etc.) has an impact on vehicle miles traveled (VMT) and public transit ridership \cite{Martin2011a, Martin2016, shaheen2018, feigon2016, lazarus2018}.}

{The emerging \gls{cav} technologies offer intriguing opportunities to enhance urban mobility and traffic safety, and the introduction of \glspl{cav} enables innovative often more responsive and efficient options for traveling which may change {the way people use }mobility services \cite{duarte2018, hancock2019}. It is likely that the wide adoption of \glspl{cav} could also affect the usage of existing infrastructure to better serve the ever-changing transportation network \cite{shaheen2017m}. While the benefits of \gls{cav} technologies on traffic flow  and safety \cite{Tientrakool2011, Shladover2012b, Fagnant2015,krueger2016,Litman2018a}, coordination in specific traffic scenarios \cite{Malikopoulos:2019a,Malikopoulos2016a,Malikopoulos2017,Malikopoulos2018a,Malikopoulos2018b,Malikopoulos2018c,Rios-Torres2015}, and energy improvement on vehicle level \cite{Malikopoulos2008b,Malikopoulos2010a,Malikopoulos2011} are well understood,} potential deployment of the \glspl{cav} for the shared mobility service has raised a number of key questions related to fleet sizing, operation strategies and the implications on mobility, urban form, and environment \cite{Shaheen2014, Zhang2015a, Rodoulis2014, Gruel2016, Ohnemus2016}.

{With the ongoing growth of shared mobility and increasing interests in \gls{cav} fleet, the convergence of emerging mobility service and technology is still evolving. Many major automakers and technology companies are launching \gls{sav} pilot projects in the US and around the world, e.g., Ford, Voyage, Waymo, Uber, and Lyft \cite{stocker2019}. While there is currently no large-scale deployment of \gls{sav} fleet, several research efforts have evaluated the impacts of the \glspl{sav}, including simulation-based evaluation on environmental impact, cost-benefit, or demand analysis, e.g., \cite{Ford2012, Burns2012, chong2013a, Fagnant2014, Bischoff2016, Dia2017, Merlin2017, Bosch2017, Moorthy2017, Metz2018, Menon2018, Foldes2018, Truong2017, tussyadiah2017, dandl2018}. There has been much contention on the potential influence of \glspl{sav} on travel behavior, urban landscape, congestion, and environment \cite{stocker2018}.} Although it seems that the required fleet size as well as the parking spaces to meet existing travel demand might drop significantly, multiple studies have indicated that full automation is likely to induce travel demand and attract new user groups, which may result in a potential increase in energy consumption, e.g., \cite{stocker2017, ross2017}. Furthermore, there have been also concerns that \glspl{sav} might attract considerable attention from public transit patrons rather than private car owners, with implications on escalating traffic congestion, if not properly managed, e.g., \cite{Chen2016b, davidson2016}.

\subsection{Scope of the Paper}

{In this paper, we review research efforts on the modeling and operations of the \gls{sav} system and try to identify potential research gaps that require further investigation. In our review, we have excluded studies on the demand estimation and travel behavior analysis of the \glspl{sav}.} {We applied the following search strings and included the papers up to date containing any combination of the keywords in the title, abstract, or keywords: 
\begin{enumerate}
    \item {Shared autonomous (electric) vehicles, shared automated vehicles, autonomous vehicle sharing;}
    \item {Autonomous carsharing, driverless carsharing, self-driving carsharing;}
    \item {Autonomous taxi, automated taxi, driverless taxi;}
    \item {Automated demand responsive transport, autonomous mobility on demand, automated    mobility on demand, autonomous mobility as a service.”}
\end{enumerate}
}

Although the exploration of benefits of \glspl{sav} is still in early stages, we note that there are many aspects in common with the conventional carsharing system (with or without the option of ridesharing). {There are several review papers providing a good summary under the umbrella of shared mobility, e.g., see \cite{Jorge2013, Agatz2012a, Furuhata2013a, Brandstatter2016, lavieri2017, jittrapirom2017, utriainen2018}. Similar review efforts on the \glspl{sav} include the study by Hao and Yamamoto \cite{Hao2018}, who focused on the features and demand aspects of the \gls{sav} system through examining the corresponding aspects of car sharing in \glspl{av}. The most recent work conducted by Stocker and Shaheen \cite{stocker2019} reviewed \gls{sav} pilots and legislation in the US, and discussed current and future development of the \gls{sav} system.}  Any such effort has obvious limitations. Space constraints limit the description of each paper in details, and thus, discussions are included only where they are important for understanding the fundamental concepts or explaining significant departures from previous work. 

\subsection{Organization of the Paper}

{The structure of the paper is organized as follows. In Section \ref{sec:mod}, we present an overview of the \gls{sav} system and modeling approaches that have been commonly adopted. We then identify major design variables and system operating parameters that are widely studied in the literature to date and summarize the research efforts in Section \ref{sec:cha}, including the problems of fleet sizing, vehicle assignment and relocation, consideration of electric vehicles, and ridesharing. In Section \ref{sec:opr}, we discuss different operation schemes of the \glspl{sav} in a mixed traffic environment that have been investigated in the literature.} Finally, we discuss research gaps and potential future research directions in Section \ref{sec:sum}.

\section{Shared autonomous vehicle system modeling}\label{sec:mod}
\glspl{sav} provide carsharing  with a way of seamlessly relocating vehicles to better match dynamic demand \cite{Fagnant2014}. {As pilot programs of \glspl{sav} are beginning to accelerate around the world, there has been an increasing interest in investigating the \gls{sav} system. In this section, we first introduce earlier work on the feasibility of statewide implementation of \glspl{sav} and system performance analysis along with the cost-benefit analysis. We then discuss two major directions in modeling and analysis of the \gls{sav} system: (1) the development of analytical models along with specific problems that include vehicle assignment and rebalancing, e.g., \cite{Zhang2014e, Zhang2016l, Spieser2014c, pavone2015}; (2) the development of agent-based models to emphasize the understanding of system performance and impact of the \gls{sav} system under different scenarios with a variety of parameters settings, e.g., \cite{Boesch2015, Chen2016a, Chen2016b, zhang2017p, Zhang2015g, Moreno2018}.}

\subsubsection{Feasibility analysis}
In an early work \cite{Ford2012}, Ford proposed a statewide \gls{sav} system in New Jersey with a grid-based network model. The author discussed different operation strategies of a \gls{sav} system at different time periods. For example, in rush hours, the \glspl{sav} would function like a personal rapid transit (PRT) system to satisfy travel demand and ease congestion, whereas during non-rush hours, the \glspl{sav} could be operated with more flexibility and provide door-to-door service. The area considered in the paper was modeled as gridded zones, where a fixed \gls{sav} station would be located at the center of each cell. Later, Brownell and Kornhauser \cite{Brownell2014b} described in detail two distinct \gls{sav} network models, i.e., PRT and the smart paratransit (SPT), and discussed the feasibility of a statewide \gls{sav} network in New Jersey. In the PRT network, fixed stations of the \gls{sav} system are established and passengers need to walk to their closest stations. Ridesharing is considered only if two passengers share the same origin-destination pair and arrive at the station within a predefined time window. {The idea behind the SPT system is that trips with close origins and/or destination will be served by one single vehicle. The vehicle moves around within the origin cell to pick up multiple passengers before traveling to the destination cell. Along the ways, the vehicle may stop at one, or more, locations to pick up or drop off passengers. In a SPT system with \glspl{av}, since the vehicle takes the place of the individual for accessing service, the distance between nodes in the transit grid could be increased.}  
Burns et al. \cite{Burns2012} conducted a cost-benefit analysis of a \gls{sav} system where the entire trip demand is satisfied by \glspl{sav}. To estimate the performance of a \gls{sav} system and compare with other systems (e.g., personal vehicle), the authors developed an analytical model with spatial queueing approach based on simplifying assumptions (e.g., uniformly distributed origins and destinations, constant trip request rate, etc). The results from three case studies showed that a \gls{sav} system is capable of providing better mobility experience at a significantly lower cost, in addition to its environmental and safety benefits. 

\subsubsection{Analytical modeling}
\begin{figure}[ht]
    \centering
    \includegraphics[width=0.45\textwidth]{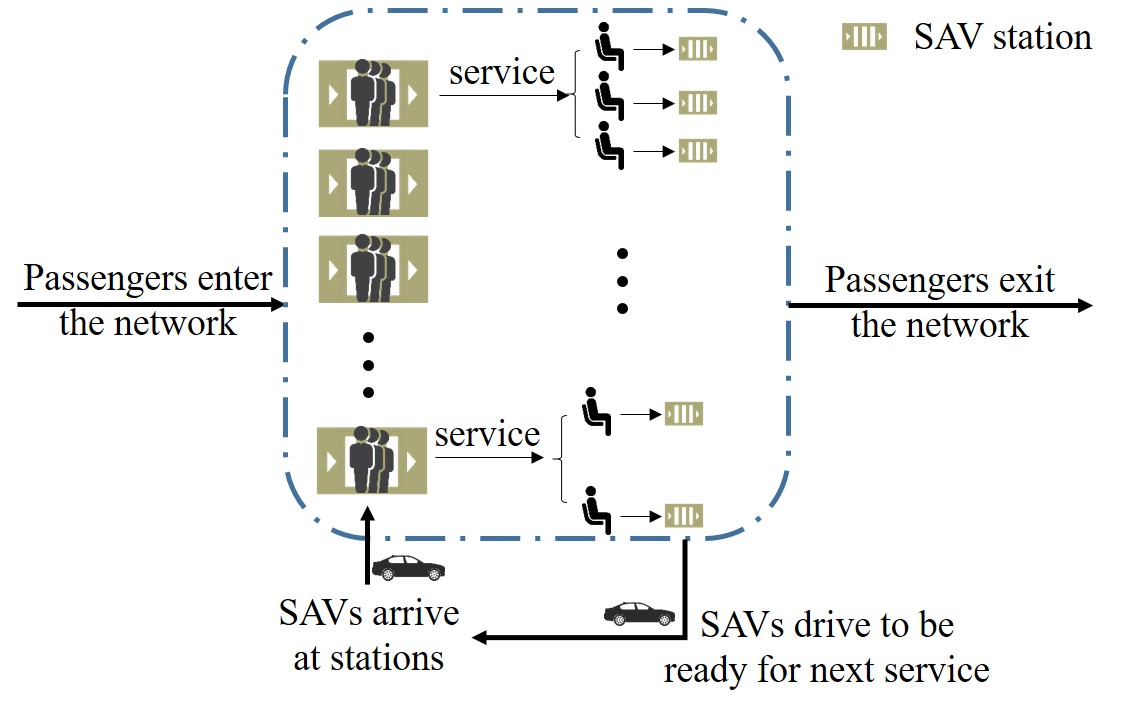}
    \caption{Shared autonomous vehicles in a queueing system.}
    \label{fig:queueing}
\end{figure}

Several research efforts reported in the literature have treated a \gls{sav} system as a spatial queueing system where passengers arrive at each station, pick up the vehicles -- if parked at the station -- and wait or leave the system, if no vehicle is available (Fig. \ref{fig:queueing}). After dropping off passengers at their destinations, vehicles either start the next service, or park, or relocate themselves to other stations, e.g., \cite{Zhang2014e, Zhang2016l, Spieser2014c, pavone2015, Rossi2018a, Iglesias2018}. For instance, Zhang et al. \cite{Zhang2015a} described a \gls{sav} network as a spatial queueing system where transportation requests queue up and are served by the \glspl{sav} in the network. The authors presented two models for \gls{sav} systems. In the first model, the authors considered a distributed approach, where the objective is to design a routing policy that minimizes the average steady-state time delay between the generation of an origin-destination pair and the time the trip is completed. In the second model, the authors considered a lumped approach -- customers are assumed to arrive at a set of stations in the network, where each customer picks up a vehicle, if available, or leaves the system, if no vehicle is parked at the station. 

\subsubsection{Agent-based modeling}
{To address the questions on the impact of \glspl{sav} on transportation mobility and investigate performance of the \gls{sav} system under various scenarios, several research efforts have also focused on developing agent-based models to evaluate the transportation network with presence of \glspl{sav} \cite{Boesch2015, Chen2016a, Chen2016b, Jager2018}.} With the advantage of modeling each individual passenger/vehicle as an agent following simple rules, complex behavior \cite{Malikopoulos2015,Malikopoulos2015b} at a macroscopic level emerges, which provides an approximation of travel behavior in the transportation systems \cite{Boesch2015}. Marczuk et al. \cite{Marczuk2015} and Azevedo et al. \cite{Azevedo2016} proposed an extension to the agent-based demand and supply model (SimMobility) for the design and evaluation of the \gls{sav} system in a multi-level simulation framework, and explored the effects of fleet size and station location for both station-based and free-floating \gls{sav} systems. Boesch and Ciari \cite{Boesch2015} discussed the advantages of MATSim (an activity-based agent-based simulation model) with the presence of \glspl{sav} and its potential applications on investigating related problems, such as the potential of \glspl{sav} complementing or competing with other transportation modes, appropriate fleet size in different transportation systems, and the demand distribution with respect to the response of different fleet sizes.

Focusing on the potential impact of a \gls{sav} system on urban parking demand, Zhang et al. \cite{zhang2017p, Zhang2015g} investigated different system operation strategies under low penetration of \glspl{sav} with an agent-based simulation model. Ridesharing and traveler's acceptance of sharing rides were also explored in the paper. The results showed a significant parking demand reduction with the \gls{sav} system -- enabling ridesharing and adding vehicle cruising options would further reduce parking demand. 
{Kondor et al. \cite{kondor2018} developed an agent-based simulation model to estimate parking demand savings with shared vehicles and \glspl{sav} for home-work commuting. Other conclusions  drawn from this study include that up to 50\% reduction in parking needs could be achieved at the expense of less than 2\% increase in VMT.}
Jager et al. \cite{Jager2018} developed an agent-based framework for a \gls{saev} system that reflect the system behavior on an operational level. Although the system has a central dispatcher, the vehicles compete for customers and make their own decisions for routing and charging. {Simulation results confirmed the feasibility of operating a \gls{sav} fleet with both high service levels and vehicle utilization. However, environmental benefits can only be expected when using renewable energy sources and enabling ride sharing features.}

\section{Shared autonomous vehicle system design variables}\label{sec:cha}
{Similar to conventional carsharing service, not only the operations of a \gls{sav} is significantly affected by the assignment and rebalancing strategies over a fleet of \glspl{sav}, mobility and environment, but also the urban landscape can be considerably impacted by the implementation strategies of a \gls{sav} system. Naturally, the problems of fleet sizing, vehicle-trip assignment, and rebalancing in a network of \glspl{sav} are the major subjects in enhancing our understanding of a \gls{sav} system, with the options of ridesharing and usage of electric vehicles that have attracted considerable attention recently. The majority of the literature to date has concentrated on how the \gls{sav} system tackles one or more of the aforementioned problems, and has aimed at enhancing our understanding about the performance and potential benefits of the network with a fleet of \glspl{sav}.} In the following subsections, we provide a summary of \gls{sav} system modeling and discuss key topics that have been investigated in previous studies regarding the \gls{sav} system.

\subsection{Fleet Sizing of a Shared Autonomous Vehicles System}
Fleet size is the major determinant of the operating cost of the \gls{sav} system. General considerations in determining the fleet size include system access, directness, sharing, and passenger waiting time \cite{Winter2016, Boesch2016}. In what follows, we summarize different approaches in addressing fleet sizing problems in a \gls{sav} system.

Fagnant et al. \cite{Fagnant2015d} simulated a \gls{sav} system in Austin area with a grid-based network model following a similar modeling framework presented in \cite{Ford2012}. In this work, a fleet of \glspl{sav} is generated in the network to ensure that passenger waiting times are within predefined bounds. A heuristic strategy is implemented to relocate vehicles such that the stock of \glspl{sav} among cells is balanced. A replacement rate of 1 \gls{sav} per 9.3 conventional vehicles was identified as appropriate for the area considered in the paper. The authors concluded that even with an excess VMT, emissions and environmental outcomes for the \glspl{sav} are still advantageous compared to those for the average US vehicle fleet.  
In the modeling framework for the \gls{sav} system developed by Winter et al. \cite{Winter2016}, the minimum fleet size and the optimal fleet size that yield the minimum system costs are determined through an iterative procedure, where the core is a simulation tool that is applied for assigning vehicles to passenger requests. Several scenarios are conducted to analyze the influence of different design parameters (e.g., vehicle capacity, operational parameters, demand level) on system performance. 

Vazifeh et al. \cite{Vazifeh2018} investigated the minimum fleet size problem of a \gls{sav} system with a network-based model. Trips based on known demand and link travel times were taken as input to construct the vehicle shareability network under the constraint of maximum trip connection time. With fully knowledge of daily trip demand, the authors found that 40\% taxis in New York City can be reduced without incurring delay to passengers, under the constraint of  15-minute maximum trip connection time. Relaxing the assumption of complete demand information, the authors concluded that if trip requests were collected at 1-minute interval, the system could be operated with a 30\% fleet reduction at a relative high level of service (i.e., above 90\% served trips within a 6-min delay).

Spieser et al. \cite{Spieser2014c} addressed two major fleet sizing problems: (1) the minimum number of vehicles needed to stabilize the workload of a \gls{sav} system and (2) the number of vehicles needed to ensure a desired level of service provided to the customers. In their paper, the \gls{sav} system is modeled as a queueing network where each region is mapped into single-server node, and each route between each pair of regions is mapped into infinite-server nodes. The vehicle rebalancing process is modeled as an arrival process of ``virtual passengers." Conducting a case study in Singapore, the paper showed that a \gls{sav} can meet the personal mobility needs of the entire population with a fleet size about one third of the total number of passenger vehicles currently in operation. 

Masoud and Jayakrishnan \cite{Masoud2017} discussed a different implementation strategy of the \gls{sav} system, with households forming clusters (i.e., neighborhoods). Each neighborhood share the ownership and ridership of a set of autonomous vehicles that serve as rental cars during their idling times. The authors focused on the optimization of the fleet size in a neighborhood and the number of rental requests for the vehicles during their idling times. Two optimization models were developed. The first model addressed the neighborhood clusters and aimed at minimizing the total number of the vehicles by considering essential trips to be satisfied for all the households in a neighborhood. The second model optimized the total number of rental requests so as to maximize extra income from idling vehicles, considering time window constraints of the owners' essential trips.

{Most of previous work has emphasized on searching for the minimum fleet size of \glspl{sav} that could provide service on the existing demand at a desired level, when replacing the existing conventional vehicle service by \glspl{sav}. We have noticed promising results from multiple papers indicating that a high replacement rate of conventional vehicles is feasible to satisfy the same level of demand. However, there is still some  work missing to assess holistically the impact of urban mobility due to potentially changing travel behavior and demand as a result of the introduction of \gls{av} in the mixed traffic environment.}

\subsection{Vehicle Assignment in a Shared Autonomous Vehicle System}

Although there is a rich body in the literature in dynamic assignment problems with various applications on taxi, paratransit, trucking services, etc, that require real-time vehicle assignment to dynamic service requests (e.g., see \cite{Spivey2004, Oncan2007, Pillac2013, laporte2009} for more details), most papers reported in the literature to date have focused on investigating \gls{sav} system performance with simplified vehicle assignment strategies (usually rule-based). 
In what follows, we present a general formulation of the vehicle assignment problem in a \gls{sav} system. Let $i \in \mathcal{M}$ be a trip request, $j \in \mathcal{N}$ be the index of a vehicle, and $x_{ij}$ equal to 1 if and only if trip $i$ is assigned to vehicle $j$, where $\mathcal{M} \subset \mathbb{N}$ is the set of trip requests and $\mathcal{N} \subset \mathbb{N}$ is the set of vehicles. The general vehicle-traveler, or vehicle-trip, assignment problem to minimize the trip assignment cost, $J_a$,  \cite{Zhang2015a, Ma2017a}:

\begin{equation}
\min J_a = \sum_i \sum_j c_{ij} x_{ij}, \\
\label{obj: assign}
\end{equation}
subject to
\begin{gather}
\sum_j  x_{ij} = 1, i \in \mathcal{M}, \label{con: traveler}\\
x_{ij} \in \{0,1\}, \forall i \in \mathcal{M}, j \in \mathcal{N},
\end{gather}
where $c_{ij}$ is the cost of assigning trip request $i$ to vehicle $j$, which could be represented by trip travel distance, travel time, or monetary cost. The trip assignment cost in \eqref{obj: assign} is evaluated at every trip assignment time step with dynamic service requests. The constraint \eqref{con: traveler} ensures that each traveler is assigned to only one vehicle. 

When assigning travelers to the nearest idling \glspl{av}, several research efforts have considered a first-come-first-served strategy, which is a heuristic approach to minimize passenger waiting time \cite{Burns2012, Zhang2015g, Fagnant2015d, Boesch2016, Shen2015}. In a paper by Fagnant and Kockelman \cite{Fagnant2014}, the \gls{sav} service area is divided into small zones, where trips are randomly generated. Every five minutes, passengers will be randomly ordered and assigned to the nearest available \gls{sav} in the same zone, up to a maximum vehicle arrival time. If such assignment fails, those passengers will be held until next assignment. 
Hyland and Mahmassani \cite{Hyland2018} investigated the underlying stochastic vehicle assignment problem for the \gls{sav} system with no shared rides. With the assumption that the fleet operator has no information of the spatial-temporal demand distribution, the authors compared different \gls{sav} assignment policies as the solution approaches to the local optimization problem at each time step. Two of the applied strategies were first-come-first-served, and the other strategies minimized traveler waiting times (under different vehicle-traveler assignment constraints). 

Hanna et al. \cite{hanna2016} examined four different methods for assigning vehicles in a \gls{sav} system: (1) a decentralized greedy matching where users are assigned to their nearest vehicles in a random order, (2) a centralized greedy matching approach ensuring that each vehicle is matched with its closest user, (3) the Hungarian minimum cost matching algorithm that minimizes passenger waiting time and unoccupied distance traveled, and (4) a minimal makespan matching algorithm which minimizes the longest distance that any vehicle must travel to a passenger. The authors showed that compared to greedy approaches, the latter two methods improved system performance through reducing unoccupied travel distance, passenger waiting time, and waiting time variation. 
 
\subsection{Vehicle Rebalancing of a Shared Autonomous Vehicle System} 
The \gls{sav} system shares similar characteristics with the carsharing system consisting of conventional vehicles \cite{Fagnant2014}. In terms of unbalanced demand distribution, both systems face the same problem of vehicle rebalancing. Two major rebalancing strategies have been investigated in the literature of carsharing with conventional vehicles including (1) operator-based vehicle relocation and (2) user-based vehicle relocation, which could potentially be adapted in addressing the same problem in the \gls{sav} system, {see \cite{Jorge2013, Jorge2014, Weikl2013, Brendel2017, santos2019}}. {However, the relocation of \glspl{sav} still have differences with that of conventional sharing vehicles, since \glspl{sav} are fully compliant and always cooperative \cite{Zhang2016g}.} Thus, due to the inherent capabilities of self-driving and self-rebalancing of a \gls{sav} system, research efforts have focused more on the problem with a centralized operator that has dispatching control over the entire \gls{sav} network, which may yield a system optimum solution for the entire system. 

We provide a general formulation to illustrate the vehicle rebalancing problem for a \gls{sav} system. Let $r_y$ be the number of idling vehicles in zone/station $y \in \mathcal{Z}$ and $r_{yz}$ be the number of rebalancing vehicles from zone/station $y$ to zone/station $z \in \mathcal{Z}$, where $\mathcal{Z} \subset \mathbb{N}$ is the total number of zones/stations in the network. Generally, the objective function $J_r$ is the total cost induced by vehicle rebalancing \cite{Wen2017, Spieser2015, Zhang2014e}:

\begin{equation}
\min J_r = \sum_y \sum_z c_{yz} r_{yz}, \\
\label{obj: relocation}
\end{equation}
subject to
\begin{gather}
\sum_z  r_{yz} = r_y, \forall y, z \in \mathcal{Z}, \label{con: trip}\\
r_{yz} \in \mathbb{N}, \forall y, z \in \mathcal{Z},
\end{gather}
where $c_{yz}$ is the cost of moving vehicles from zone/station $y$ to zone/station $z$, which could be represented by trip travel distance, travel time, or monetary cost.
In a system with dynamic trip requests, \eqref{obj: relocation} will be evaluated at every rebalancing time step and \eqref{con: trip} defines the total rebalancing vehicles from zone/station $y$ should equal the number of idling vehicles in the zone. 

Targeting at the problem of unbalancing demand and supply, Pavone et al. \cite{Pavone2012} addressed the vehicle relocation problem for a mobility-on-demand system, optimizing the rebalancing assignment that minimizes the number of vehicles to be moved. Using a fluid model of the system, the authors showed that the optimal rebalancing policy can be found as the solution to a linear program, under which every station reaches an equilibrium where there are excess vehicles and no waiting customers. Based on this study, Zhang and Pavone \cite{Zhang2016l} presented a queueing-theoretical approach and provided the solution to an offline optimal rebalancing problem. Later, Wen et al. \cite{Wen2017} extended the research by incorporating door-to-door service and ridesharing option in a free-floating \gls{sav} system. 
From the fleet operator's perspective, Spieser et al. \cite{Spieser2015} investigated the vehicle rebalancing problem in a \gls{sav} system by quantifying the operation cost as a function of fleet size, demand loss and utilization rate, and analyzed the impact of fleet size on demand loss, vehicle utilization rate, and vehicle rebalancing miles traveled. {H{\"o}rl et al. \cite{horl2019} evaluated performance of four heuristic and optimal rebalancing policies for a \gls{sav} system in an agent-based simulation environment, and suggested that the utilization of intelligent demand forecasts and rebalancing algorithms would be crucial for a \gls{sav} system to be competitive with private vehicles.}

Through simulation based evaluation, recent work focused on the impact of vehicle rebalancing strategies in a \gls{sav} system. Zhu and Kornhauser \cite{Zhu2017} investigated the rebalancing strategies for the \gls{sav} system in New Jersey and their effects on the fleet size and level of service provided in scenarios where all non-walking travel demand is served by \glspl{sav}. Shared trips are served by vehicles of different capacities (i.e., 3, 6, 15, and 50 passengers). Two rebalancing strategies are developed based on known demand. In the first approach, vehicles are moved at the end of the day to make sure that there are enough vehicles at each station that satisfy the demand at the beginning of the day. In the second approach, vehicles are relocated as needed to fill in any station without enough vehicles. The authors also evaluated the performance of the statewide \gls{sav} system with varying fleet sizes, in terms of passenger waiting time and rebalancing trip lengths. The results showed that one \gls{sav} could possibly replace more than six traditional vehicles while the demand could still be well served. 

{Fagnant and Kockelman \cite{Fagnant2014} investigated the operation of \glspl{sav} through an agent-based model and focused on the implications of travel and environmental impacts of \glspl{sav} under a mixed traffic condition. Addressing the imbalanced demand patterns, the authors proposed several relocation strategies to balance vehicle supply and reduce future traveler wait times: (1) relocating vehicles based on expected demand and (2) relocating vehicles to balance stock based on predicted supply.} Marczuk et al. \cite{Marczuk2016} developed a simulation framework for rebalancing an one-way \glspl{sav} system in SimMobility environment. The proposed fleet management center is responsible for passenger-to-vehicle assignment, vehicle routing and re-balancing. Three vehicle relocation strategies were proposed for the system: (1) no rebalancing as the baseline scenario, (2) offline rebalancing that minimizes the number of rebalancing trips, and (3) online rebalancing that minimizes the total time/effort spent for rebalancing per rebalancing interval. Winter et al. \cite{Winter2017} analyzed the impacts of different relocation strategies of a \gls{sav} system in a simulated generic grid network. Five vehicle relocation strategies were tested, including remaining idle, random shuffling, returning to original location, moving based on demand anticipation, and moving to balance vehicle stock over the network. {In the simulation framework, the fleet size of the \gls{sav} system is given as an input, and vehicles are dispatched through a rule-based strategy. Performance measures such as average passenger utility, average waiting time, and the ratio of vehicle driving time were examined. The simulation showed that remaining idle strategy would be the most efficient in terms of passenger waiting time, yet the worst performer considering link occupancy and parking turnover rates. In contrast, strategies aimed at distributing vehicles yielded higher parking turnover rates but showed lower service efficiency. } 
{In light of these results, the authors extended the study by imposing the constraints of limited parking facilities in the evaluation of the above five heuristic relocation strategies for idle \glspl{sav}, and examined the potential impact of \glspl{sav} on urban traffic in terms of congestion, parking consumption and mode shift \cite{winter2018}.}

{As discussed in the above papers, e.g., \cite{Marczuk2016, Zhu2017, Winter2017, winter2018}, depending on the objectives and targeting performance measures, the rebalancing strategy to be applied in a \gls{sav} system may be different. The operation of a fleet of \gls{sav} is considerably affected by the applied relocation strategy or a combination of strategies, considering the inter-dependencies among parking demand, traffic condition, and user mode choice. Although current research efforts emphasize rebalancing strategies in an isolated \gls{sav} system, the externalities should be analyzed in more depth to enhance the understanding of traffic dynamics with the implementation of \gls{sav} service.}

\subsection{The Usage of Electric Vehicles in a Shared Autonomous Vehicle System}

A significant amount of work has focused on the use of electric vehicles in a \gls{sav} system to achieve larger energy and emission savings for a greener transportation system \cite{Boesch2016, Chen2016a, Chen2016b}. Considering the range of electric vehicles, there is a number of constraints in a \gls{saev} system. For instance, a vehicle may need to visit a charging station after dropping off passengers. There may be instances that vehicles have to turn down trip requests and drive to charging stations instead, resulting in different vehicle-trip assignment strategies \cite{Kang2016, Loeb2018, Loeb2019}. 

Based on the work in \cite{Zhang2015a}, Zhang et al. \cite{Zhang2016g} presented a model predictive control (MPC) approach to optimize vehicle scheduling and routing in a \gls{saev} system, considering vehicle charging constraints. Compared to other control algorithms of a \gls{sav} system (i.e., nearest-neighbor dispatch, collaborative dispatch, Markov redistribution, real-time rebalancing), the authors concluded with a case study in New York City that the MPC algorithms outperformed the other strategies in terms of average customer waiting times.  

Chen et al. \cite{Chen2016a, Chen2016b} addressed the operations of a \glspl{saev} with an agent-based model based on the work reported in \cite{Fagnant2014} and \cite{Fagnant2015d}. The emphasis of this research is the performance analysis of a fleet of SAEVs under various vehicle range and charging infrastructure scenarios. The authors also explored the pricing schemes of a SAEV system when competing against other modes (i.e., private human-driven vehicles and city bus service), and found that with higher SAEV penetration rate, the private vehicle replacement rate by the SAEVs increases, leading to improved system performance. {Similarly, the study by Bauer et al. \cite{bauer2018} predicted battery range and charging infrastructure requirements of a fleet of \glspl{saev} operating on Manhattan island with an agent-based model. The authors also conducted sensitivity analysis of the cost and the environmental impact of providing \gls{saev} service with a wide range of changes in cost components (e.g., battery type, vehicle type, etc.). The study indicated that instead of battery range, the major challenge to introducing \glspl{saev} may be building sufficient charging infrastructure.}

Kang et al. \cite{Kang2016} developed a framework for a \gls{saev} system that consists of demand forecasting, fleet assignment, electric vehicle designing, and charging station locating modules. The fleet assignment module determines the optimal vehicle assignment and charging schedules, and the charging station locating module decides the optimal charging station locations. The system-level objective is to maximize  service profit for the operator, through optimizing decision variables including fleet size, number of charging stations, electric powertrain design, membership fee, and vehicle rental fee. The locations of charging stations are selected with a p-median model from a pool of predetermined candidates. A comparison between a \gls{sav} system and a \gls{saev} system was conducted in terms of cost and benefit under different scenarios (e.g., varying gas prices and charging station installation costs), showing that a \gls{saev} system would be more profitable for most of the scenarios. Although both systems are marketable, the optimized \glspl{saev} required longer waiting times than optimized \glspl{sav} due to the constraints of vehicle range and charging issues. 

Iacobucci et al. \cite{Iacobucci2018} developed a simulation model to evaluate a \gls{saev} system interacting with passengers and charging at designated stations based on a heuristic charging strategy. The potential utilization of the \gls{saev} system as an operating reserve provider and its performance in response to grid operator requests were evaluated. The authors concluded that the proposed system could reduce the required fleet size as compared to private vehicles while providing a comparable level of transportation service with low break-even prices. {Later, based on the work presented in \cite{Zhang2016g}, the authors developed a framework for the optimization of charging scheduling and vehicle routing and relocation for a fleet of \glspl{saev} \cite{iacobucci2019}. The proposed framework consists of two layers of optimization model: over longer time scales, the charging scheduling {optimization minimizes waiting times and electricity costs,} while over shorter time scales, vehicle routing and relocation are optimized under charging constraints. The authors reported that a substantial reduction in charging costs was yielded from the proposed framework without significantly affecting passenger waiting times, as well as the potential of \glspl{saev} to offer energy storage to the grid and avoid grid congestion.}

{In summary, the introduction of electric vehicles in the \gls{sav} system offers a large potential to further enhance environmental benefits. However, constraints such as vehicle range and charging facility locations add more dynamics into the system, and multiple studies suggested that the infrastructure and charging scheduling are the key influencing factors of system performance of a fleet of \glspl{saev}. Considerably work has focused on the performance analysis of \gls{saev} system as compared to the \gls{sav} system, through evaluating the impact of vehicle range, charging infrastructure, as well as electricity costs \cite{Chen2016a, Chen2016b, bauer2018}. Considering charging constraints, several research efforts have also emphasized on re-examining vehicle routing and relocation strategies as well as optimizing charging locations \cite{Zhang2016g, Kang2016, miao2019}. Recently, the option of vehicle-to-grid as well as the integrated planning of power grid and shared mobility service has also attracted considerable attention \cite{Iacobucci2018, iacobucci2019}, to improve the perception of \glspl{saev} and ensure sustainable commutes within the notion of smart cities \cite{Oldenbroek2017}.}

\subsection{The Option of ridesharing in a Shared Autonomous Vehicle System}

\begin{figure}[htb]
\centering

\begin{subfigure}[b]{0.45\textwidth}
\includegraphics[width=\textwidth]{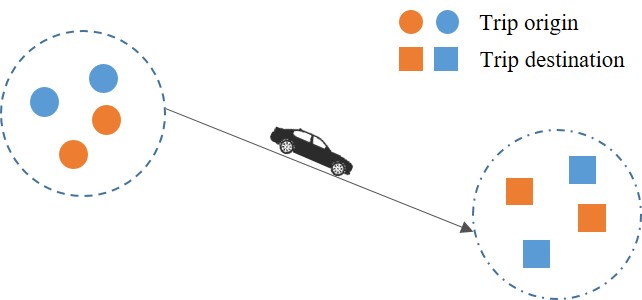}
\caption{}
\label{fig:r_a}
\end{subfigure}%

\begin{subfigure}[b]{0.45\textwidth}
\includegraphics[width=\textwidth]{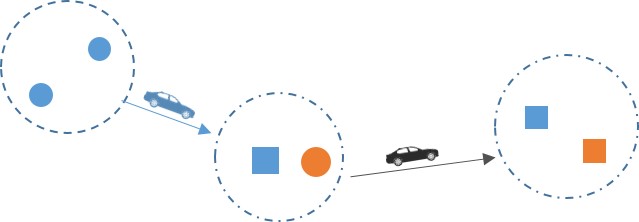}
\caption{}
\label{fig:r_b}
\end{subfigure}%

\caption{ridesharing in the shared autonomous vehicle system: a) trip combination; b) trip chaining.}
\label{fig:ridesharing}
\end{figure}

The problems of ridesharing and carsharing are usually decoupled in the existing literature \cite{Samaranayake2017}. Recently research efforts started exploring the option of ridesharing in a \gls{sav} system, e.g., \cite{Levin2017,Hyland2018a, Fagnant2018, tsao2019}. By allowing ridesharing, the fleet size may be further reduced to provide a desired level of service to the passengers, although the total VMT probably might increase \cite{Burghout2015, Gurumurthy2018}. There are generally two types of ridesharing as illustrated in Fig.~\ref{fig:ridesharing}: (a) trip combining neighboring origins and destinations (Fig.~\ref{fig:r_a}) and (b) trip chaining based on trip temporal and spatial characteristics (Fig.~\ref{fig:r_b}). {We consider here ridesharing as the option of serving multiple passengers in a single vehicle trip, or trip chain, in the \gls{sav} system, and emphasize the impact of opening up ridesharing options in the \gls{sav} service, without detailing the operation modes and strategies for ridesharing. } Considering different system objectives (e.g., minimizing total VMT, minimizing total travel time, or maximizing served trips) and various system constraints (e.g., time window and seat constraints), there has been work on the \gls{sav} system with the option of ridesharing and the evaluation of different ridesharing strategies against network performance.

\begin{table*}[t] 
\centering
\caption{Approaches in Shared Autonomous Vehicle System Modeling.}
\begin{tabular*}{.9\textwidth}{ccc}
Approach & Topic & Reference  \\
\midrule
Optimization & Fleet sizing & \cite{Spieser2014c,Vazifeh2018,Zhang2015a,Ma2017a,Iglesias2017,Masoud2017,Kang2016} \\
& Vehicle routing / trip assignment& \cite{Levin2017a,Liang2017,Rossi2018a,Liang2018,Zhang2016g, Ma2017a,Hyland2018a,Liang2016, liang2018o, de2016, iacobucci2019, Samaranayake2017, tsao2019}\\
& Vehicle rebalancing / relocation & \cite{Spieser2015,Wen2017,Zhang2016l,Rossi2018a,Iglesias2017} \\
& Other considerations & \cite{Kang2016,Liang2016, miao2019} \\
\midrule
Simulation Evaluation&Fleet sizing&\cite{Bischoff2016,Brownell2014b,Burns2012,Burghout2015,Boesch2016,Dia2017,Ford2012,Winter2016,Marczuk2015}\\
&Vehicle routing / trip assignment&\cite{Burghout2015,Shen2015,Rigole2014,Azevedo2016,Pinto2018,Jager2018}\\
&Vehicle rebalancing / relocation&\cite{Brendel2017,Burghout2015,Fiedler2018,Winter2017,Dia2017,Rigole2014,Scheltes2017,Fagnant2014,Fagnant2015d,Marczuk2016,zhang2017p,winter2018, horl2019}\\
&Ridesharing&\cite{Levin2017,Farhan2018,Fagnant2018,Zhang2015f,Zachariah2014,Lu2018,shen2016a,shen2018,zhang2017p,chen2019, Gurumurthy2018}\\
&Pricing scheme&\cite{Horl2016,Liu2017,Chen2016b,shen2016a, dandl2018}\\
&Transit integration / mode choice&\cite{Horl2016,Martinez2017,Liu2017,Azevedo2016,Scheltes2017,Winter2016,Shen2016b,shen2018,Zachariah2014,Chen2016b,Pinto2018, snelder2019, wen2018, basu2018}\\
&Electric vehicles&\cite{Brendel2017,Farhan2018,Chen2016a,Iacobucci2018,Scheltes2017,Loeb2018,Jager2018,bauer2018, dandl2018, Loeb2019}\\
\midrule
\end{tabular*}
\label{tab:summary}
\end{table*}

Levin et al. \cite{Levin2017} analyzed the possibility of ridesharing in a \gls{sav} system where passengers could select the first arrived vehicle regardless of occupancy. The authors found that \glspl{sav} with the choice of ridesharing may cause more congestion due to additional miles traveled for detouring. Zhang et al. \cite{Zhang2015f, Zhang2015g} applied an agent-based model to evaluate the performance and potential benefits of a \gls{sav} system with dynamic ridesharing. In a grid-based simulation network, a centralized operator monitors real-time trip requests and \gls{sav} status as well as manages trip assignment for the \gls{sav} system, where ridesharing option is evaluated against passenger's willingness and travel cost. Their work suggested that dynamic ridesharing in a \gls{sav} system could potentially lead to reduced vehicle ownership, parking demand, and emissions.

Hyland and Mahmassani \cite{Hyland2018a} compared the performance of a \gls{sav} system with and without ridesharing option in terms of the ability to handle demand surges. In this paper, the mathematical formulations of the vehicle assignment with/without ridesharing were presented and the solutions were derived with a rolling-horizon approach. The simulation results indicated that the \gls{sav} with ridesharing service improved system performance in response to demand surges.

Based on the vehicle rebalancing strategies tested in \cite{Fagnant2014}, Fagnant and Kockelman \cite{Fagnant2018} considered the option of dynamic ridesharing in a simulated \gls{sav} system. With the case study of a 24-mile by 12-mile region in Austin, the authors concluded that dynamic ridesharing in a \gls{sav} system was able to limit excess VMT from the \gls{sav} system, reduce passenger waiting times (under the constraint that ridesharing should not increase travel time of current passengers by more than 40\%), and yield an enhanced level of service. 

{Farhan and Chen \cite{Farhan2018} discussed the impacts of ridesharing on the operational efficiency of \glspl{saev} with a discrete-time simulation model. Both the fleet size and number of charging stations are determined during simulation. In their research, the travelers are grouped into clusters based on spatial criteria, and the ride-share matching problem is formulated as a vehicle routing problem minimizing system-wide vehicle miles traveled under time window constraint. The results indicated that allowing a second passenger in ridesharing yielded marginal benefit of fleet size and charging station reduction. Although more passengers in shared trips reduced the required fleet size and number of charge stations, passenger waiting times increased due to ridesharing (i.e., reduced level of service).} 

\section{Shared Autonomous Vehicle System Operation} \label{sec:opr}

{Although the majority of the literature has been focused on examining the feasibility and performance of the \gls{sav} service as an isolated system, there is an increasing interest towards the investigation of more realistic operational scenarios} related to the \glspl{sav}. Recent research efforts have also focused on answering questions such as: ``How will the \gls{sav} system {perform} in a mixed traffic environment?"  ``What will be the mobility impact {of} integrating the \glspl{sav} with other modes of transport?" In this section, we focus on different operational aspects of a \gls{sav} system, and summarize the studies that consider realistic and mixed traffic conditions.

\subsection{Operation in a Realistic Traffic Environment}

The majority of the aforementioned work has addressed the \gls{sav} system with full \gls{sav} penetration or without considering background traffic. Only a few papers have focused on the  impact of congestion of \glspl{sav}, e.g., \cite{Levin2017,Liang2017,Rossi2018a,lamotte2017,maciejewski2017}. For example, to investigate the impact of \glspl{sav} on mobility, Levin et al. \cite{Levin2017} presented a general event-based framework for simulating the operations of a \gls{sav} system with existing traffic models. Considering 100\% penetration of \glspl{sav}, the authors found that under certain scenarios (e.g., with the option of dynamic ridesharing), a smaller fleet of \glspl{sav} performed better than a larger fleet due to lower congestion in the network. Maciejewski and Bischoff \cite{maciejewski2017} evaluated the impact of a city-wide introduction of \glspl{sav} on traffic congestion through an agent-based simulation model, focusing on the analysis of traffic congestion under different \gls{sav} penetration rates. Under the assumption of increased road capacity due to \gls{av} operations, their work showed that despite increased traffic volume, a fleet of \gls{sav} could have a positive effect on traffic at a penetration rate as low as 20\%. 

Levin \cite{Levin2017a} developed a linear programming formulation for vehicle routing problem in the \gls{sav} system, where traffic flow was modeled through the link transmission model. The results showed that asymmetric demand (e.g., demand during peak periods) could lead to significantly rebalancing trips and greater congestion than uniformly distributed demand pattern. Since more vehicles might cause additional congestion on roadway network, it is important for the \gls{sav} system to plan for different traffic patterns.
Liang et al. \cite{Liang2017} proposed an integer programming model to define the routing of the \glspl{sav} based on profit maximization function, where travel times on the links varied with the flow of \glspl{sav} (without any background traffic). Later in \cite{Liang2018}, the authors applied the algorithm for trip assignment and dynamic routing in the city of Delft, the Netherlands with a rolling horizon scheme. Assuming that the operator of a \gls{sav} fleet has the choice of accepting or rejecting trip requests according to a profit maximization function, the analysis showed that taking into account the impact of dynamic travel time led to different results of satisfied trips and VMT, and ultimately affected overall operator profit and network congestion level.

Rossi et al. \cite{Rossi2018a} studied the routing and rebalancing problem of \glspl{sav} in congested transportation networks, where a \gls{sav} system is modeled in a network flow framework such that vehicles are represented as flows in a road network. The objective of the routing problem is to minimize the weighted sum of passenger trip travel times and vehicle rebalancing travel times considering network capacity. The objective of the rebalancing problem is to optimize rebalancing paths such that traffic congestion is minimized. Through numerical studies on real-world traffic data, the authors showed that the proposed real-time routing and rebalancing algorithm yielded lower customer waiting time by avoiding excess congestion on the road, compared to point-to-point rebalancing algorithms where no underlying road network is assumed. 

{Through an agent-based model, Fagnant and Kockelman \cite{Fagnant2014} investigated the operation of \glspl{sav} and focused on the implications of travel and environmental impacts of \glspl{sav} under a mixed traffic condition. Instead of 100\% penetration of \glspl{sav}, the authors considered the transportation system with a small market share of \glspl{sav} (i.e., around 3.5\%). The simulation results under different scenarios (e.g., varying trip generation rates, network congestion levels, \gls{sav} fleet size, etc) indicated that each \gls{sav} can substitute around eleven conventional vehicles at the cost of 10\% more VMT, and the overall emissions savings are expected to be sizable for most emission species.} 

\subsection{Operating in a Multi-Modal Environment}

Based on the discussion in the previous sections, it seems clear that \glspl{sav}, compared to personal owned human-driven vehicles, have significant advantages for individuals as well as for the transportation system in terms of mobility, safety, and energy savings (especially with \glspl{saev}), {e.g., \cite{Marczuk2015,Chen2016a,Jager2018, Fagnant2014,Bischoff2016, Dia2017}. A combination of \glspl{sav} with other transportation modes such as public transportation, however, might impose different conclusions \cite{Martin2016, shaheen2018, feigon2016, lazarus2018}. Although \glspl{sav} could be utilized in the way to facilitate the first and last mile transport \cite{Malikopoulos2019b} and promote the use of public transportation system (e.g., \cite{boersma2018, basu2018}), \glspl{sav} may also divert passengers away from transit systems due to their capability of providing door-to-door services (e.g., \cite{wadud2016, basu2018}). }

\subsubsection{Shared autonomous vehicles as a complement of public transit}
Early research efforts have explored the performance of integrating the \gls{sav} system with transit systems. For example, based on the same network in New Jersey as in \cite{Ford2012}, Zachariah et al. \cite{Zachariah2014} simulated a system of \glspl{sav} where the train network is preserved and treated as an integral part of the system. 
Using \glspl{sav} as a complementary service of a train system, Liang et al. \cite{Liang2016} presented an optimization model to define the service area of a \gls{sav} system for first/last mile transport that maximizes the profit of the \gls{sav} operator. {Later in \cite{liang2018o}, the authors designed a \gls{sav} system providing shuttle service between a major train station and city area, considering the competition between \glspl{sav} or other modes (e.g., biking or walking), as well as the impact of traffic congestion on mode split. With the objective of minimizing total travel time, the authors developed an optimization model to decide the best fleet size and price rate for the \gls{sav} system.}

Shen et al. \cite{Shen2016b, shen2018} explored the feasibility of integrating \glspl{sav} in the public transportation system to improve the first/last mile connectivity. With a simplified simulation model without considering traffic congestion where the demand for the \gls{sav} system was assumed to be 10\% of the original bus demand, the study showed that by enabling ridesharing, the integrated service was able to reduce average passenger travel time and ease traffic through less occupancy of road resources.
Scheltes and de Almeida Correia \cite{Scheltes2017} studied the \gls{saev} system providing last-mile service for a train line. In the simulation model, vehicle assignment in response to traveler request followed a first-come-first-served model. The scenarios of short-term pre-booking, vehicle relocating, and opportunity charging were also explored. The results showed that compared to bicycle and walking as last mile transportation modes, the \gls{saev} system was able to reduce average passenger travel time and waiting time, especially when pre-booking option was enabled. 

{Wen et al. \cite{wen2018} proposed a systematic approach to design and simulate an integrated system of \glspl{sav} and public transit. The authors emphasized {that} the \gls{sav} operation is designed to be transit-oriented with the purpose of supporting existing public transit service. In an agent-based simulation platform, the interaction between service operator and travelers {is modeled} with a set of system dynamics equations, such that the decisions of both parties could be captured in the system. The authors suggested that encouraging ridesharing, allowing in-advance requests, and combining fare with transit would be useful to enable service integration and promote sustainable travel.}
Pinto et al. \cite{Pinto2018} proposed a simulation framework integrating a travel mode choice model and a dynamic transit assignment model to assess the impacts of a suburban first-mile \gls{sav} system on transit demand. Similarly, Martinez and Viegas \cite{Martinez2017} presented an agent-based model to evaluate the impact of the \glspl{sav} in the city of Lisbon, Portugal. In their simulation model, current travel demand is served by two types of \glspl{av} that compete with each other, i.e., a \gls{sav} providing door-to-door service with the choice of ridesharing and an autonomous minibus that replaces current bus service without any transfers for users. The simulation results revealed positive mobility impact of \glspl{sav} especially when introducing the autonomous minibus into the network. 

\subsubsection{Shared autonomous vehicles as a competitor of public transit}
Liu et al. \cite{Liu2017} simulated transportation patterns in Austin network with a system of \gls{sav} from a mode-choice perspective. A user-equilibrium based dynamic traffic assignment model was applied in  simulation environment. The study focused on travelers' mode choices with the presence of \glspl{sav}. In a mixed traffic environment, where private human-driven vehicles, public transit, and \glspl{sav} coexist, the study analyzed the impacts of the \gls{sav} system on energy consumption and emissions under different \gls{sav} penetration rates and \gls{sav} rental fees. {Based on the sensitivity analysis of rental fees, the authors found that if the \gls{sav} fare rate is low enough, \gls{sav} users might travel more than private vehicle users. Therefore, although the use of \glspl{av} {is expected to result} in energy savings and emission reduction, the extra VMT by \glspl{sav} could compromise such environmental benefits. The mode choice results indicated that, for travelers who do not own a private vehicle, \glspl{sav} are preferable for short-distance trips compared to public transit. However, demand shifting from public transit would be a concern once the \glspl{sav} become available in the study area.} 
H{\"{o}}rl \cite{Horl2016} conducted a similar study and investigated the \gls{sav} service in a multi-modal traffic simulation environment. The simulation results in the test scenario raised the following two concerns: (1) the introduction of \glspl{sav} led to increased VMT and, moreover, (2) \glspl{sav} attracted public transportation users rather than private car owners.

{Snelder et al. \cite{snelder2019} developed a simulation framework to assess both direct and indirect impacts of \glspl{av} and \glspl{sav} in a mixed traffic environment. To capture demand elasticities, the network fundamental diagram was combined with mode choice models. Furthermore, {the spatial impact was} also modeled as an exogenous input to the framework via a percentage of relocated inhabitants per lane use type. The simulation results showed that a shift to \glspl{sav} could be expected. However, the improved accessibility for many residents could result in a significant increase in vehicle trips (and also in VMT), which might impose negative effects on traffic condition. Similar conclusions were drawn from the study on the effects of full automation with the possibility of trip chaining of household trips, yet in a scenario where most vehicles are still privately owned \cite{de2016}.}

{In summary, findings of multiple studies indicate that although the introduction of \glspl{sav} in the transportation system might improve mobility and safety, it could result in enormous changes of travel behavior, mode choice, car ownership, and possibly transportation infrastructure and urban form. A holistic assess of the impact of the {\gls{sav} systems} on urban mobility and related social implications might be challenging at the moment as \glspl{sav} are still evolving. However, \gls{sav} service could possibly have negative impact on traffic congestion and be strongly competitive with public transit without appropriate incentive mechanisms.}

\section{Outlook and Future Directions} \label{sec:sum}
\subsection{Concluding Remarks}
In this paper, we summarized current research efforts in \gls{sav} systems that have been reported in the literature to date. Although the \gls{sav} system have many aspects in common with the conventional carsharing system, the inherent characteristics of self-driving and self-rebalancing with \glspl{sav} further enhance free-floating carsharing service and increase the stochasticity of the system internally. {Externally, the introduction of \glspl{av} in the transportation network could change fundamentally traffic patterns in the future. The complexity of traffic and urban dynamics, thus, places considerable uncertainty in terms of both short-term and long-term impacts of the system \cite{milakis2017}.}

The majority of research efforts has considered a system either of full \gls{sav} penetration rate or without any traffic, and compared its performance with the conventional mobility systems  (in terms of fleet size requirement, energy implications, VMT, passenger travel times, etc). Among these research efforts, agent-based modeling is one of the major approaches to evaluate network performance of a \gls{sav} system and assess potential impacts of the system. {Several research efforts have focused on developing optimization models to address the following questions: (1) ``what is the minimum fleet size to provide a desired level of service?" (2) ``What is the optimal vehicle assignment strategy to minimum passenger travel time?" (3) ``What is the optimal vehicle relocation strategy to minimize the number of rebalancing trips without inducing waiting delay?" In general, the \gls{sav} system could benefit from the cooperative characteristics of the fleet -- the connectivity and automation embedded in the system open up the opportunities for a central controller to apply optimal operation strategies to achieve global optimum against different network design objectives.} 

Although previous research has aimed at enhancing our understanding of the \gls{sav} systems, there are still open issues to be addressed. For example, most papers consider the \gls{sav} system with fixed stations whereas free floating \gls{sav} systems have not been thoroughly investigated. Within a \gls{sav} system, the optimal fleet sizing problem to maintain a minimum required level of service or to ensure a desired level of service is still under-explored. The considerations of different vehicle assignment and relocation strategies, or the option of ridesharing further increase the complexity of the problem. So far most papers have applied heuristics for the implementation of \glspl{sav} to solve these problems and focused more on assessing potential benefits of a \gls{sav} system. 

\subsection{Future Research}
There are several directions for future research considering the gaps in the work reported in the literature to date. Although previous work has addressed the replacement ratio of \glspl{sav} to conventional private vehicles, the majority of the results are derived with existing demand patterns in an isolated system. The problem of modeling the \gls{sav} system with presence of other transportation modes, as either a complement or competing mode, needs further investigation. Especially, the following questions still remain unanswered: (1) ``What is the network performance of a \gls{sav} system in a realistic transportation network?" (2) ``How much improvement in the level of service in a transportation network can be achieved with an integrated \gls{sav} system?"  To address these challenges, it is necessary to study the operational strategies (e.g., optimal fleet size/vehicle assignment/relocation strategy, etc) which would yield the minimum and/or desired level of service of the transportation network. Furthermore, in an environment where massive amount of data could be collected from vehicles and infrastructure, what we used to model as uncertainty become an additional input. With the advent of information and communication technologies, better utilizing available information for optimal operational strategies requires novel solutions to reduce dimensions and {to overcome issues associated with data in high-dimensional spaces.}

With all possible mobility service options enabled by CAVs, one particular question that still remains unanswered is ``how demand pattern or travel behavior will eventually change?" With the shared mobility choices (and enhanced convenience with \glspl{sav}), there is already an evidence of an increase of induced demand (e.g., more night travels, or trips shifted from transit demand). However, little research has been conducted on investigating the impact of the emerging \gls{sav} system on the vulnerable population, while a systematic framework of providing accessibility to a variety of social groups is still missing. Meanwhile, the nature of self-driving and self-rebalancing of a \gls{sav} system also implies potential changes on land use. For example, the implications of a \gls{sav} system on urban parking spaces is still under-explored. Thus, the long-term impact of shared mobility system on urban transportation systems is still an open question.


\ifCLASSOPTIONcaptionsoff
  \newpage
\fi



\bibliographystyle{IEEEtran}
\bibliography{SharedMobilityRef}%


\begin{IEEEbiography}
[{\includegraphics[width=1.1in,height=1.25in,clip,keepaspectratio]{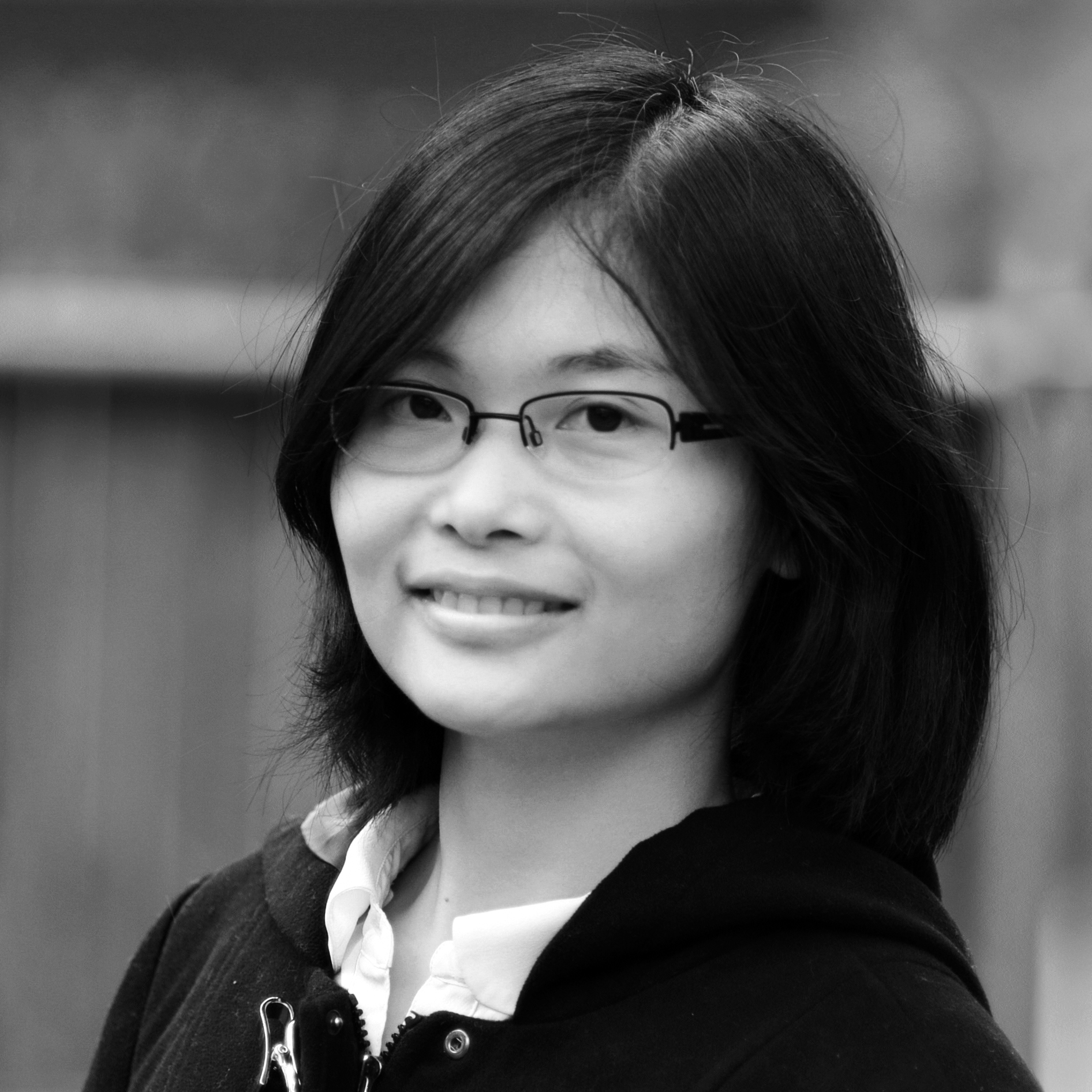}}]
{Liuhui Zhao} (M2017) received the B.S. degree in Resources Science and Technology from Beijing Normal University, Beijing, China, in 2009, the M.S. degree from Department of Geography at the University of Alabama in 2011, and the Ph.D. degree in Transportation Engineering from New Jersey Institute of Technology in 2016. She is currently a Postdoctoral Researcher in the Information and Decision Science (IDS) Laboratory at the University of Delaware leading research projects on emerging transportation systems. She has participated in various research projects on connected automated vehicles, intelligent transportation systems, traffic and transit operations. Her research interests lie within the areas of intelligent transportation systems, shared mobility, and connected automated vehicles.
\end{IEEEbiography}

\begin{IEEEbiography}[{\includegraphics[width=1.1in,height=1.25in,clip,keepaspectratio]{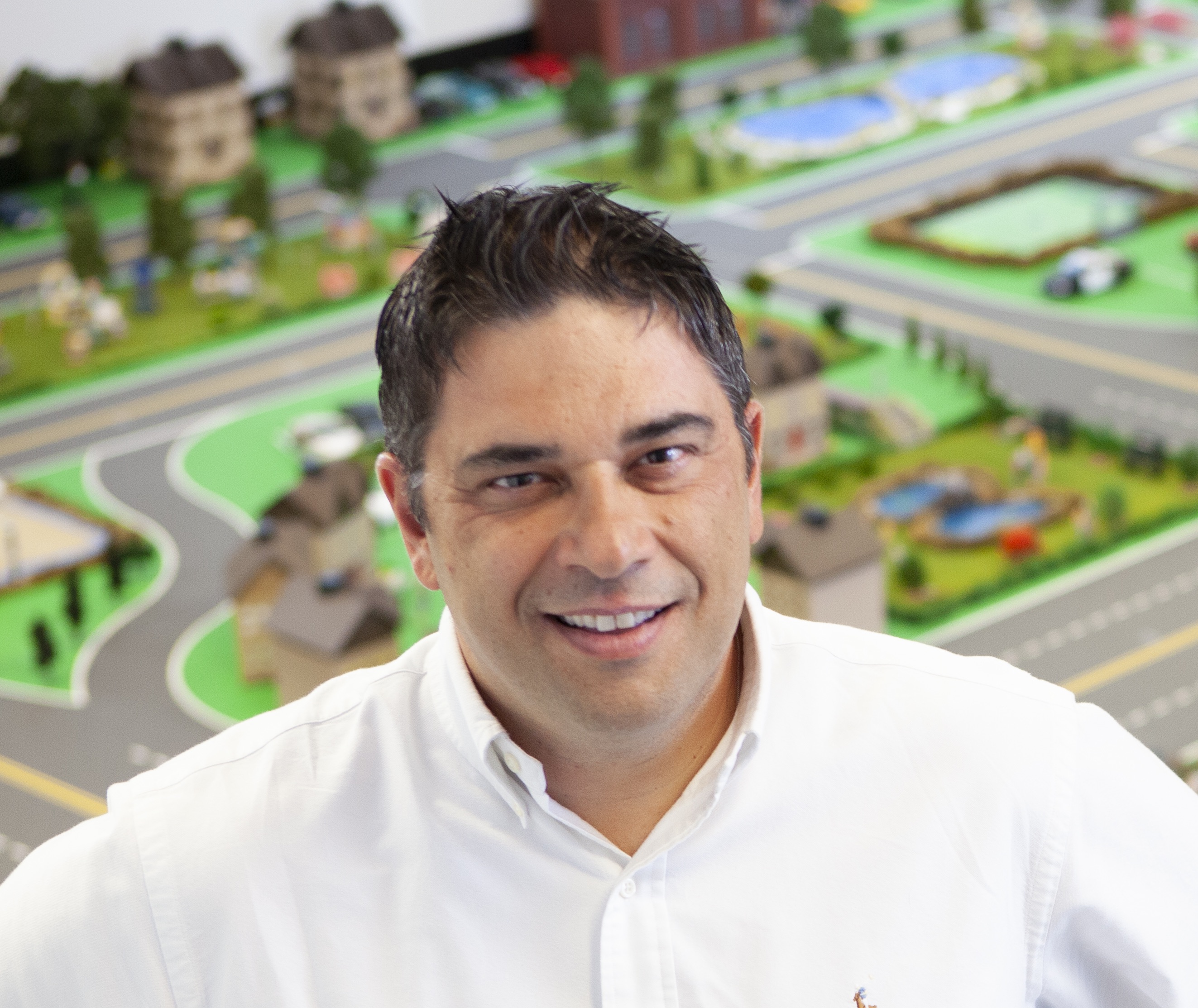}}]{Andreas A. Malikopoulos}
(M2006, SM2017) received a Diploma in Mechanical Engineering from the National Technical University of Athens, Greece, in 2000. He received M.S. and Ph.D. degrees from the Department of Mechanical Engineering at the University of Michigan, Ann Arbor, Michigan, USA, in 2004 and 2008, respectively. He is the Terri Connor Kelly and John Kelly Career Development Associate Professor in the Department of Mechanical Engineering at the University of Delaware (UD) and the Director of the Information and Decision Science (IDS) Laboratory. Before he joined UD, he was the Deputy Director and the Lead of the Sustainable Mobility Theme of the Urban Dynamics Institute at Oak Ridge National Laboratory, and a Senior Researcher with General Motors Global Research \& Development. His research spans several fields, including analysis, optimization, and control of cyber-physical systems; decentralized systems; and stochastic scheduling and resource allocation problems. The emphasis is on applications related to sociotechnical systems, energy efficient mobility systems, and sustainable systems. He is currently an Associate Editor of the IEEE Transactions on Intelligent Vehicles and IEEE Transactions on Intelligent Transportation Systems. He is a member of SIAM, AAAS, and a Fellow of the ASME.
\end{IEEEbiography}



\end{document}